# Closures in $\aleph_0$-categorical bilinear maps

Andreas Baudisch

Humboldt-Universität zu Berlin
Institut für Mathematik
D-10099 Berlin
e-mail: baudisch@mathematik.hu-berlin.de

**Abstract**

Alternating bilinear maps with few relations allow to define a combinatorial closure similarly as in [2]. For the $\aleph_0$-categorical case we show that this closure is part of the algebraic closure.

## 1 Introduction

A.H. Lachlan asked in [6] whether every $\aleph_0$-categorical stable theory is $\omega$-stable. 1989 E. Hrushovski [5] refuted this conjecture. He constructed an $\aleph_0$-categorical stable graph that is not $\omega$-stable.

$\aleph_0$-categorical stable groups have been investigated by W. Baur, G. Cherlin, and A. Macintyre [3] and by U. Felgner [4]. They proved independently that $\aleph_0$-categorical stable groups are nilpotent by finite. Furthermore in [3] is proved, that every $\aleph_0$-categorical, $\omega$-stable group is abelian by finite. Therefore W. Baur, G. Cherlin, and A. Macintyre formulated the following conjecture weaker than Lachlan's conjecture:

**BCM-Conjecture**  An $\aleph_0$-categorical stable group is abelian by finite.

This problem is still open. For a solution it is sufficient to show that every $\aleph_0$-categorical stable alternating bilinear map is trivial by finite.
In the next section we introduce the notion of an alternating bilinear map with few relations. Then we can apply the ideas of [2]. It is possible to define a closure as in [2], which is not a pregeometry in every case. Under the assumption of $\aleph_0$-categoricity we show that this closure is part of the algebraic closure.

## 2 Alternating bilinear maps

Let $p$ be a prime. We define $\mathcal{B}(p)$ to be the category of the following structures $\langle M, W, \beta \rangle$ where $M$ and $W$ are vector spaces over the field $I\!F_p$ with $p$ elements and $\beta$ is an alternating bilinear map of $M \times M$ into $W$, where the image of $\beta$ generates $W$. The morphisms of $\mathcal{B}$ are the monomorphisms between these structures. If $p$ is fixed we write only $\mathcal{B}$. If we speak about bilinear maps we mean the structures in some $\mathcal{B}(p)$. Since $\aleph_0$-categorical stable groups are nilpotent by finite (see above) and of finite exponent it is easy to show the following

**Lemma 2.1** *The BCM-Conjecture is true if and only if for every $p$ the $\aleph_0$-categorical stable bilinear maps in $\mathcal{B}(p)$ are trivial by finite.*

If $p$ is fixed, then in [2] a category $\mathcal{S}$ is introduced equivalent to $\mathcal{B}$. If $M$ is any $\mathbb{F}_p$ vector space, then there is a "free bilinear" map $\langle M, \Lambda^2 M, \wedge \rangle$ for $M$ in $\mathcal{B}$ such that for all $\langle M, W, \beta \rangle$ in $\mathcal{B}$:

$$\begin{array}{ccc} M \times M & \xrightarrow{\beta} & W \\ & \searrow_{\wedge} & \uparrow f_\beta \\ & & \Lambda^2 M \end{array}$$

where $f_\beta$ is a vector space homomorphism. The vector space $\Lambda^2 M$ (equipped with the associated bilinear map $\wedge$) is called the exterior square of $M$. The kernel of $f_\beta$ on $\Lambda^2 M$ is denoted by $N(\beta)$. The elements of the auxiliary category $\mathcal{S}$ are the pairs $(M, N(\beta))$. If we work in a fixed bilinear map $\langle M, W, \beta \rangle$ and $H$ is a subspace of $M$, then we denote by $(H, N(H))$ the pair corresponding to the restriction of $\beta$ to $H \times H$. Then $N(H) = N(\beta) \cap \Lambda^2 H$, where $\Lambda^2 H$ is considered as a subspace of $\Lambda^2 M$ in the canonical way.

**Definition** $\langle M, W, \beta \rangle$ *has few relations*, if there is a natural number $k$ such that $\dim(N(H)) \leq k \dim(H)$ for all finite subspaces $H$ of $M$. dim denotes the linear dimension of the vector spaces.

Note that $\dim(\Lambda^2 H) = \dim(H)(\dim(H) - 1)/2$. The BCM-Conjecture implies:

**Conjecture** There is no $\aleph_0$-categorical stable bilinear map with few relations in any $\mathcal{B}(p)$.

To investigate bilinear maps with few relations we introduce a combinatorial closure similarly as in [2].

## 3 The combinatorial closure

Let $\langle M, W, \beta \rangle$ be a bilinear map in some $\mathcal{B}(p)$ with few relations. Then there is some $k$ such that $\dim(H) \leq k \dim(N(H))$ for every finite subspace $H$ of $M$.

**Definition** For every finite subspace $H$ of $M$ we define
$$\delta_k(H) = \dim(H) - \frac{1}{k}\dim(N(H)).$$

Then we have $\delta_k(H) \geq 0$. For the following development we use the corresponding proofs in [2].

**Lemma 3.1** *Let $H$ and $K$ be finite suspaces of $M$.*

i) $\delta_k(H+K) \leq \delta_k(H) + \delta_k(K) - \delta_k(H \cap K)$.

ii) *If $\delta_k(H+K) > \delta_k(K)$, then $\delta_k(H \cap K) < \delta_k(H)$.*

iii) *If $\delta_k(H+K) \geq \delta_k(K)$, then $\delta_k(H \cap K) \leq \delta_k(H)$.*

**Definition** Let $H$ be a finite subspace of $M$. Then $H$ is *selfsufficient*, if for every finite subspace $K \supseteq H$ we have $\delta_k(K) \geq \delta_k(H)$.

**Lemma 3.2** *If $H$ and $K$ are selfsufficient, then $H \cap K$ is selfsufficient.*

Using this lemma and $\delta_k(H) \geq 0$ for all finite $H \subseteq M$ we can define:

**Definition**

i) For every finite subspace $H$ the selfsufficient closure $CSS(H)$ of $H$ is the intersection of all finite selfsufficient subspaces that contain $H$.

ii) Assume $A$ and $B$ are finite subsets of $M$. We say $A \subseteq \mathrm{cl}_k(B)$, if $\delta_k(CSS(\langle B \rangle)) = \delta_k(CSS(\langle A \cup B \rangle))$.

iii) $d_k(A) = \delta_k(CSS(\langle A \rangle))$.

Note that $CSS(H)$ exists, is finite, and definable. $CSS(\langle A \rangle)$ defines a closure operator on finite sets. Hence for an infinite subset $B$ of $M$ we can define:

**Definition**

i) $CSS(\langle B \rangle)$ is the union of all $CSS(\langle A \rangle)$ where $A$ is a finite subset of $B$.

ii) $B$ is selfsufficient, if $B$ is a subspace and for every finite subset $A$ of $B$ there exists a finite selfsufficient subspace $B_0$ of $B$ such that $A \subseteq B_0 \subseteq B$.

**Lemma 3.3** $\mathrm{cl}_k$ *is a closure operator on finite sets.*

**Proof.** i) $A \subseteq \mathrm{cl}_k(A)$.
This is clear by definition.

ii) $A \subseteq B$ implies $\mathrm{cl}_k(A) \subseteq \mathrm{cl}_k(B)$.
As mentioned above $A \subseteq B$ implies
$$CSS(\langle A \rangle) \subseteq CSS(\langle B \rangle).$$
Assume $C \subseteq \mathrm{cl}_k(A)$. Then
$$\delta_k(CSS(\langle B \cup C \rangle)) = \delta_k(CSS(\langle B \cup (C \cup A) \rangle)) \leq \delta_k(CSS(\langle B \rangle) \cup CSS(\langle C \cup A \rangle))$$
since by definition
$$\delta_k(CSS(\langle B \cup (C \cup A) \rangle)) \leq \delta_k(K)$$
for every finite subspace $K$ with $\langle B \cup (C \cup A) \rangle \subseteq K$. Using this and Lemma 3.1 we obtain
$$\begin{aligned}\delta_k(CSS(\langle B \cup C \rangle)) &\leq \delta_k(CSS(\langle B \rangle)) + \delta_k(CSS(\langle C \cup A \rangle)) \\ &\quad - \delta_k(CSS(\langle B \rangle) \cap CSS(\langle C \cup A \rangle)).\end{aligned}$$
By Lemma 3.2
$$CSS(\langle B \rangle) \cap CSS(\langle C \cup A \rangle)$$
is a selfsufficient subspace between the selfsufficient subspaces $CSS(\langle A \rangle)$ and $CSS(\langle C \cup A \rangle)$. Since by assumption $\delta_k(CSS(\langle A \rangle)) = \delta_k(CSS(\langle A \cup C \rangle))$ we get
$$\delta_k(CSS(\langle B \rangle) \cap CSS(\langle C \cup A \rangle)) = \delta_k(CSS(\langle C \cup A \rangle)).$$
Hence
$$\delta_k(CSS(\langle C \cup B \rangle)) \leq \delta_k(CSS(\langle B \rangle)).$$
By selfsufficiency of $CSS(\langle B \rangle)$ we have equality as desired.

iii) $\mathrm{cl}_k(\mathrm{cl}_k(A)) = \mathrm{cl}_k(A)$.
By ii) and i) $\mathrm{cl}_k(\mathrm{cl}_k(A)) \supseteq \mathrm{cl}_k(A)$. Assume $B \subseteq \mathrm{cl}_k(A)$ and $C \subseteq \mathrm{cl}_k(B)$. We have to show $C \subseteq \mathrm{cl}_k(A)$. We have
$$\delta_k(CSS(\langle B \cup A \rangle)) = \delta_k(CSS(\langle A \rangle)).$$
By ii) $C \subseteq \mathrm{cl}_k(B)$ implies $C \subseteq \mathrm{cl}_k(B \cup A)$. We get
$$\delta_k(CSS(\langle C \cup B \cup A \rangle)) = \delta_k(CSS(\langle B \cup A \rangle)).$$
Hence
$$\delta_k(CSS(\langle C \cup B \cup A \rangle)) = \delta_k(CSS(\langle A \rangle)),$$
as desired. □

By Lemma 3.3 we can define for infinite subsets $B$:

**Definition**

i) Let $A$ be a finite subset of $M$. Then $A \subseteq \text{cl}_k(B)$, if there exists a finite subset $B_0$ of $B$ such that $A \subseteq \text{cl}_k(B_0)$.

ii) Let $H$ and $K$ be subspaces of $M$ such that $H+K/H$ is finite and $H$ is selfsufficient. Then

$$\delta_k(K/H) = \min\{\delta_k(K_1) - \delta_k(K_1 \cap H) : \\ K_1 \cap H \text{ is selfsufficient, } K_1 \text{ finite, } K_1 + H = K + H\}.$$

Since we have assumed that $H$ is selfsufficient the definition makes sense. By selfsufficiency of $H$ there is a finite $K_0$ such that $K_0 + H = K + H$, $K_0 \cap H$ is selfsufficient, and $N(K_0 + H) = N(K_0) + N(H)$. Then $\delta_k(K/H) = \delta_k(K_0) - \delta_k(K_0 \cap H) = \dim(K + H/H) - \frac{1}{k}\dim(N(K + H)/N(H))$.

## 4 The free case

Let $\langle M, \Lambda^2 M, \wedge \rangle$ be the free alternating bilinear map for the vector space $M$. We show the well-known result that this structure is not $\aleph_0$-categorical. This will be done since parts of the proof are used later. The basis of the investigation of $\langle M, \Lambda^2 M, \wedge \rangle$ is the following:

**Fact** If $\{a_i : i < \kappa\}$ is a basis of $M$, then $\{a_i \wedge a_j : i < j < \kappa\}$ is a basis of $\Lambda^2 M$.

We call the $a_i \wedge a_j$'s the basic commutators over $A$.

**Lemma 4.1** Let $A = \{a_i : i < 2m\}$ be a set of linearly independent elements of $M$. Let $w = \sum_{i<m}(a_{2i} \wedge a_{2i+1})$. Let $E = \{e_i : i < n\}$ be a set of linearly independent elements in $\langle A \rangle$ with $n < m$. Then $w \notin \Lambda^2\langle E \rangle$.

**Proof.** Let $s$ be $2m - n$. Then there exists a subset $I$ of $\{0, \ldots, 2m-1\}$ of size $s$ such that $A_0 = \{a_i : i \in I\}$ is linearly independent over $\langle E \rangle$ and $A_0 \cup E$ is a basis of $\langle A \rangle$. Then there are a injective function $f$ from $\{0, \ldots, n-1\}$ onto $\{0, \ldots, 2m-1\} \setminus I$ and elements $r^i, r^i_j \in \mathbb{F}_p$ such that

$$(*) \qquad e_i = r_i a_{f(i)} - \sum_{j \in I} r^i_j a_j \text{ and } r_i \neq 0 \text{ for } i < n.$$

Since $n < m$ there is some $\ell$ such that $2\ell$ and $2\ell + 1$ are elements of $I$. W.l.o.g. we assume $\ell = 0$. Using a basis transformation of $E$ we can furthermore suppose that there are at most one $i_0$ such that $r^{i_0}_0 \neq 0$ and at most one $i_1$ such that $r^{i_1}_1 \neq 0$. Now assume $w \in \Lambda^2\langle E \rangle$. We show a contradiction.

By $(*)$ we can rewrite the presentation of $w$ over $E$ as a presentation of $w$ over $A$. Therefore $r^{i_0}_0 \neq 0$, $r^{i_1}_1 \neq 0$, and $i_0 \neq i_1$ is necessary. But then the obtained presentation of $w$ over $A$ contains the basic commutators $a_0 \wedge a_{f(i_1)}$ and $a_1 \wedge a_{f(i_0)}$ a contradiction.□

**Corollary 4.2** *Let $M$ be an infinite vector space over $\mathbb{F}_p$. Then the elementary theory of $\langle M, \Lambda^2 M, \wedge \rangle$ is not $\aleph_0$-categorical.*

**Proof.** Let $\{a_i : i < \kappa\}$ be a basis of $M$. By assumption $\omega < \kappa$. Let $g(0) = 1$ and $g(i+1) = 2g(i) + 1$. We define
$$w_{g(i)} = \sum_{j<g(i)} (a_{2j} \wedge a_{2j+1}).$$
By Lemma 4.1 the elements $w_{g(i)}$ determine different orbits of our structure. □

## 5 Minimal extensions

**Definition** Let $H$ and $K$ be selfsufficient subspaces of $M$ such that $H \subseteq K$. We call $K$ a minimal extension of $H$ if $K/H$ is finite, $\delta_k(K/H) = 0$, and for every subspace $L$ with $H \neq L \neq K$ and $H \subseteq L \subseteq K$ we have $\delta_k(L/H) > 0$.

In this section we prove two lemmas that are essential for the final results.

**Lemma 5.1** *Assume $K$ and $L$ are selfsufficient extensions of $H$ and $H$ is also selfsufficient. If $K$ is a minimal extension of $H$, then $K \subseteq L$ or $K$ and $L$ are linearly independent over $H$.*

**Proof.** Assume neither $K \subseteq L$ nor $K \cap L = H$. Let $n$ be $\dim(K/H)$ and $K_1 = K \cap L$. Then $H \subseteq K_1 \subseteq K$ and $H \neq K_1 \neq K$. Since $K$ is a minimal extension of $H$ there is a basis $a_0 \ldots, a_{r-1}$ of $K$ over $K_1$ and a basis $\Phi_0, \ldots, \Phi_{s-1}$ of $N(K)$ over $N(K_1)$ such that $0 < r < n$ and $s > k \cdot r$. But then $a_0 \ldots a_{r-1}$ is a basis of $K + L$ over $L$ and $\Phi_0 \ldots, \Phi_{s-1}$ can be extended to a basis of $N(K+L)$ over $N(L)$. Hence $\delta_k(K+L/L) < 0$ a contradiction to the selfsufficiency of $L$. □

**Lemma 5.2** *Let $H$ be selfsufficient and let $K$ and $L$ be selfsufficient extensions of $H$ such that $K, L \subset \mathrm{cl}_k(H)$ and $L \cap K = H$. Then $N(L+K) = N(L) + N(K)$.*

**Proof.** It is sufficient to show the assertion for the case where $H$, $L$, and $K$ are finite. By the assumption
$$\delta_k(H) = \delta_k(K) = \delta_k(L) \leq \delta_k(K+L).$$
By Lemma 3.1
$$\delta_k(L+K) \leq \delta_k(L) + \delta_k(K) - \delta_k(L \cap K).$$
Since $L \cap K = H$ and $\delta_k(H) = \delta_k(K)$
$$\delta_k(L+K) \leq \delta_k(L) = \delta_k(H).$$
Since $H$ is selfsufficient we have
$$\delta_k(L+K) = \delta_k(L) = \delta_k(K) = \delta(H).$$
Therefore
$$N(L+K) = N(L) + N(K). \qquad \square$$

# 6 The closure is finite

We prove the result of the paper

**Theorem 6.1** *Let $\langle M, W, \beta \rangle$ be an $\aleph_0$-categorical bilinear map in $\mathcal{B}$ with few relations. If $H$ is a definable selfsufficient subspace of $M$, then $\mathrm{cl}_k(H)/H$ is finite. Hence $\mathrm{cl}_k(H)$ is definable with the parameters of $H$ and $\mathrm{cl}_k(H) \subseteq \mathrm{acl}(H)$.*

**Corollary 6.2** *Let $\langle M, W, \beta \rangle$ be an $\aleph_0$-categorical bilinear map in $\mathcal{B}$ with few relations. If $H$ is a finite subspace of $M$, then $\mathrm{cl}_k(H)$ is finite and part of $\mathrm{acl}(H)$.*

**Proof.** $CSS(H)$ is finite and $H$-definable. Therefore we can apply Theorem 6.1. $\square$

**Proof of Theorem 6.1.**

First step:

(1) Let $H$ be any definable selfsufficient subspace of $M$. If $\mathrm{cl}_k(H) \setminus H \neq \emptyset$, then there is a minimal extension of $H$.

**Proof of (1).** We use the fact that (1) is true for finite selfsufficient $H$. Assume $a \in \mathrm{cl}_k(H) \setminus H$. Then there is some finite selfsufficient subspace $H_0$ of $H$ such that $a \in \mathrm{cl}_k(H_0)$. $H \cap \mathrm{cl}_k(H_0)$ is selfsufficient by Lemma 3.2. Hence w.l.o.g. $H_0 = \mathrm{cl}_k(H_0) \cap H$. $H_0$ has a minimal extension $K \subseteq \mathrm{cl}_k(H_0)$. Then $H + K$ is a minimal extension of $H$. $\square$

(2) Assume for every definable selfsufficient $H$ in $M$ that $H$ has only finitely many minimal extensions. Then $|\mathrm{cl}_k(H)/H|$ is finite.

**Proof of (2).** If $\mathrm{cl}_k(H) \setminus H \neq \emptyset$, then $H$ has a minimal extension. Let $H_1$ be the subspace of $\mathrm{cl}_k(H)$ generated by the finitely many minimal extensions of $H$. Then $|H_1/H|$ is finite. By Lemma 5.1 $H_1$ is generated by a finite number of minimal extensions of $H$ that are linearly independent over $H$. By Lemma 5.2 $H_1$ is again selfsufficient and $\delta_k(H_1/H) = 0$. $H_1$ is definable using the parameters of $H$ only. If $H_1 = \mathrm{cl}_k(H)$, then we are done. Otherwise by (1) there is a minimal extension of $H_1$ in $\mathrm{cl}_k(H) = \mathrm{cl}_k(H_1)$. We obtain $H_2$ generated by the finitely many minimal extension of $H_1$. As above $|H_2/H|$ is finite, $H_2$ is selfsufficient, $\delta_k(H_2/H) = 0$, and definable over the parameters of $H$. As long as $H_{i-1} \neq \mathrm{cl}_k(H)$ we can construct $H_i$ such that $H_{i-1} \subseteq H_i \subseteq \mathrm{cl}_k(H)$, $H_{i-1} \neq H_i$, $|H_i/H|$ finite, $H_i$ selfsufficient, $\delta_k(H_i/H) = 0$, and $H_i$ definable over the parameters of $H$.
By $\aleph_0$-categoricity there is some $i$ such that $H_i = \mathrm{cl}_k(H)$. It follows $|\mathrm{cl}_k(H)/H|$ is finite as desired. $\square$

By (2) it is sufficient to show that every selfsufficient definable subspace $H$ has only a finite number of minimal extensions. We assume that this is not true for $H$ and show a contradiction.

Every minimal extension of $H$ can be written as $H + K$ where $K$ is finite, $H \cap K$ is selfsufficient, $N(H + K) = N(H) + N(K)$ and $K$ is a minimal extension of $H \cap K$. We call $K$ a finite realization of a minimal extension of $H$. By $\aleph_0$-categoricity finitely many isomorphism types $(K, K \cap H)$ of finite realizations suffice to describe all possible minimal extension of $H$. Hence by our assumption there infinitely many different finite realizations $K_0\, K_1\, K_2 \ldots$ of minimal extensions of $H$ such that all $(K_i, K_i \cap H)$ are of the same isomorphism type. By Lemma 5.1 we can furthermore assume that the $K_i$'s are linearly independent over $H$.

For each $m$ we choose $a_i \in K_i \setminus H$ and define
$$w_m = \sum_{i<m} \beta(a_{2i}, a_{2i+1}).$$

Then $w_m$ fulfils a formula $\varphi_m(x)$ over the parameters of $H$ that says:

$$\exists X_0 \ldots X_{2m-1} \, \exists x_0 \ldots x_{2m-1}$$
$$\Big( \bigwedge_{i<2m} \text{"}X_i \text{ has the fixed isomorphism type over } H\text{"}$$
$$\wedge \text{"the } X_i\text{'s are linearly independent over } H\text{"}$$
$$\wedge \bigwedge_{i<2m} x_i \in X_i \setminus H$$
$$\wedge\, x = \sum_{i<m} \beta(x_{2i}, x_{2i+1}) \Big).$$

We get the desired contradiction if we show that $w_m$ does not fulfil $\varphi_n(x)$ for $2n < m$. Assume $2n < m$ and $\varphi_n(w_m)$ is true. Then there are finite realizations $L_i$ $(i < 2n)$ of minimal extensions of $H$ of the fixed isomorphism type over $H$ such that $L_0, \ldots, L_{2n-1}$ are linearly independent over $H$ and there are $b_i \in L_i \setminus H$ with
$$w_m = \sum_{i<n} \beta(b_{2i}, b_{2i+1}).$$

Let

$$\begin{aligned} A &= \{a_0, \ldots, a_{2m-1}\}, \\ B &= \{b_0, \ldots, b_{2n-1}\}, \\ U &= \Big(\sum_{i<2m} K_i\Big) + H, \text{ and} \\ V &= \Big(\sum_{i<2n} L_i\Big) + H. \end{aligned}$$

By Lemma 5.2 $U$ and $V$ are selfsufficient. Using Lemma 5.1 there exists $K_{i_0}, \ldots, K_{i_{s-1}}$ linearly independent over $V$ such that
$$K_i \subseteq \Big(\sum_{t<s} K_{i_t}\Big) + V \text{ for } i < 2m.$$

Then by Lemma 5.1
$$s \geq 2m - 2n.$$
Furthermore $U + V = \sum_{t<s} K_{i_t} + V$. Let $I$ be $\{i_t : t < s\}$. Now we choose a basis $D$ of $V + U$ such that

i) For $j \in I$  $D$ contains a basis $D_j$ for $K_j$ over $H$ that contains $a_j$.

ii) $D = \bigcup_{j \in I} D_j \cup E \cup F$ where $F$ is a basis of $H$ and $E$ is a basis of $V$ over $H$.

iii) For every $\ell \notin I$ there are some $r_j^\ell$ for all $j \in I$, $d_\ell \in \langle \bigcup_{j \in I}(D_j \setminus \{a_j\})\rangle$, and $c_\ell$ such that $c_\ell = a_\ell - \sum_{j \in I} r_j^\ell a_j$ and $c_\ell - d_\ell \in E$.

Let $A_0 = \{a_i : i \in I\}$ and $C_0 = \{c_i : i \notin I\}$. By i) and iii)  $A_0 \cup C_0$ is a basis of $\langle A \rangle$. Now we consider the situation in $\langle M, \Lambda^2 M, \wedge \rangle$. For elements in $\Lambda^2 M$ we use the same notation as for their preimages in $\langle M, W, \beta \rangle$. Since $|C_0| \leq 2n < m$ we can apply Lemma 4.1 and get
$$w \notin \Lambda^2 \langle C_0 \rangle.$$
Hence there is a presentation of $w$ as a linear combination of basic commutators over $A_0 \cup C_0$ of the following form

(3) $\qquad w = r(a \wedge e) + u$

where $a \in A_0$, $r \neq 0$, $e \in (A_0 \setminus \{a\}) \cup C_0$ and $u \in \Lambda^2 \langle (A_0 \setminus \{a\}) \cup C_0 \rangle$. W.l.o.g. $a = a_0$. To get a representation of $w$ with respect to $D$ we have to replace $c_\ell$ by $c_\ell - d_\ell + d_\ell$ where $c_\ell - d_\ell$ and $d_\ell \in D$ ($\ell \notin I$). We use the following identities:

$$\begin{aligned} a_j \wedge c_\ell &= a_j \wedge (c_\ell - d_\ell) + (a_j \wedge d_\ell) \\ c_\ell \wedge c_k &= (c_\ell - d_\ell) \wedge (c_k - d_k) + (c_\ell - d_\ell) \wedge d_k + d_\ell \wedge (c_k - d_k) + d_\ell \wedge d_k. \end{aligned}$$

With respect to $D$ (3) has the form

(4) $\qquad w = r(a_0 \wedge (f_1 + f_2)) + u^*$

where $r \neq 0$, $f_1 \neq 0$, $f_1 \in \langle D \setminus (D_0 \cup F)\rangle$, $f_2 \in \langle D_0 \setminus \{a_0\}\rangle$, and $u^* \in \Lambda^2 \langle D \setminus \{a_0\}\rangle$.
Besides $w$ we consider $v = \sum_{i<n} b_{2i} \wedge b_{2i+1}$ in $\langle M, \Lambda^2 M, \wedge \rangle$.
Then $v \in \Lambda^2 \langle E \cup F \rangle$.
Now we go back to $\langle M, W, \beta \rangle$. In $\langle M, W, \beta \rangle$ we have $w = v$. This means $w - v = r(a_0 \wedge (f_1 + f_2)) + u^* - v \in N(D)$. By Lemma 5.2
$$N(\langle D \rangle) = N(K_0 + H) + N(\langle D \setminus D_0 \rangle).$$
Since $u^* - v \in \Lambda^2 \langle D \setminus \{a_0\}\rangle$, $a_0 \wedge f_2 \in \Lambda^2(K_0 + H)$, and $f_1 \notin K_0 + H$ we see that $a_0 \wedge f_1$ cannot be cancelled. It is impossible that
$$r(a_0 \wedge (f_1 + f_2)) + u^* - v \in N(K_0 + H) + N(\langle D \setminus D_0 \rangle),$$
as desired. $\qquad\square$